\documentclass[11pt]{article}

\usepackage{textcomp}
\usepackage{ccfonts} 
\usepackage{euler}
\usepackage{amssymb}
\usepackage{amsmath}
\usepackage{amscd}
\usepackage{latexsym}
\usepackage{graphicx}   
\usepackage[all]{xy}    
\usepackage{psfrag} 
\usepackage{color}
\usepackage{epsfig}

\definecolor{trustcolor}{rgb}{0.99,0.0,0.0}

\newtheorem{theorem}{Theorem}

\newtheorem{definition}[theorem]{Definition}

\newtheorem{lemma}[theorem]{Lemma}

\newtheorem{proposition}[theorem]{Proposition}

\newtheorem{rema}{Remark}
\newenvironment{remark}{\begin{rema} \rm}{\end{rema}}
\newtheorem{exam}{Example}

\def\E{\mathbb E}
\def\EXP{{\E}}

\def\PROB{{\mathbb P}}

\def\IND#1{{\mathbb I}_{{\left[ #1 \right]}}}
\def\I{{\mathbb I}}

\def\argmin{\mathop{\rm arg\, min}}

\def\sgn{\mathop{\rm sgn}}

\def\blackslug{\hbox{\hskip 1pt \vrule width 4pt height 8pt depth 1.5pt
\hskip 1pt}}
\def\qed{\quad\blackslug\lower 8.5pt\null\par\medskip}

\newcommand{\X}{{\cal X}}

\def\S{{\cal S}}

\def\K{{\cal K}}

\def\Var{{\rm Var}}
\newcommand{\RR}{\mathbb{R}}

\newcommand{\Prob}[1]{\mathbb{P}\left\{ #1 \right\} }

\allowdisplaybreaks

\begin{document}
\title{Ranking the Best Instances}
\author{St\'ephan Cl\'{e}men\c{c}on \\
MODALX - Universit\'{e} Paris X Nanterre\\
\&\\
Laboratoire de Probabilit\'{e}s et Mod\`{e}les Al\'{e}atoires\\
UMR CNRS 7599 - Universit\'{e}s Paris VI et Paris VII\\
\and Nicolas Vayatis\footnote{\noindent Address of corresponding
author: Nicolas Vayatis, Laboratoire de Probabilit\'es et
Mod\`eles Al\'eatoires - Universit\'e Paris 6 - 175, rue du
Chevaleret - 75013 Paris, France - Email: {\tt
vayatis@ccr.jussieu.fr}} \\
Laboratoire de Probabilit\'{e}s et Mod\`{e}les Al\'{e}atoires\\
UMR CNRS 7599 - Universit\'{e}s Paris VI et Paris VII}

\maketitle

\begin{abstract}
We formulate the local ranking problem in the framework of
\textit{bipartite ranking} where the goal is to focus on the best
instances. We propose a methodology based on the construction of
real-valued scoring functions. We study empirical risk
minimization of dedicated statistics which involve empirical
quantiles of the scores. We first state the problem of {\em
finding} the best instances which can be cast as a classification
problem with mass constraint. Next, we develop special performance
measures for the local ranking problem which extend the Area Under
an ROC Curve (AUC/AROC) criterion and describe the optimal
elements of these new criteria. We also highlight the fact that
the goal of ranking the best instances cannot be achieved in a
stage-wise manner where first, the best instances would be
tentatively identified and then a standard AUC criterion could be
applied. Eventually, we state preliminary statistical results for
the local ranking problem.

\noindent {\bf Keywords:} Ranking, ROC curve and AUC, empirical
risk minimization, fast rates.

\noindent {\bf Running title:} Ranking the best instances
\end{abstract}

\pagebreak

\section{Introduction}
The first takes all the glory, the second takes nothing. In
applications where ranking is at stake, people often focus on the
best instances. When scanning the results from a query on a search
engine, we rarely go beyond the one or two first pages on the
screen. In the different context of credit risk screening, credit
establishments elaborate scoring rules as reliability indicators
and their main concern is to identify risky prospects especially
among the top scores. In medical diagnosis, test scores indicate
the odds for a patient to be healthy given a series of
measurements (age, blood pressure, ...). There again a particular
attention is given to the "best" instances not to miss a possible
diseased patient among the highest scores. These various
situations can be formulated in the setup of bipartite ranking
where one observes i.i.d. copies of a random pair $(X,Y)$ with $X$
being an observation vector describing the instance (web page,
debtor, patient) and $Y$ a binary label assigning to one
population or the other (relevant vs. non relevant, good vs. bad,
healthy vs. diseased). In this problem, the goal is to rank the
instances instead of simply classifying them. There is a growing
literature on the ranking problem in the field of Machine Learning
but most of it considers the Area under the ROC Curve (also known
as the AUC or AROC) criterion as a measure of performance of the
ranking rule \cite{CorMoh04, FISS03, RCMS05, AGHHPR05}. In a
previous work, we have mentioned that the bipartite ranking
problem under the AUC criterion could be interpreted as a
classification problem with pairs of observations \cite{CLV05}.
But the limit of this approach is that it weights uniformly the
pairs of items which are badly ranked. Therefore it does not
permit to distinguish between ranking rules making the same number
of mistakes but in very different parts of the ROC curve. The AUC
is indeed a global criterion which does not allow to concentrate
on the "best" instances. Special performance measures, such as the
Discounted Cumulative Gain (DCG) criterion, have been introduced
by practitioners in order to weight instances according to their
rank \cite{JarKek00} (see also \cite{Rud06, CosZha06}) but
providing theory for such criteria and developing empirical risk
minimization strategies still is a very open issue. In the present
paper, we extend the results of our previous work in \cite{CLV05}
and set theoretical grounds for the problem of local ranking. The
methodology we propose is based on the selection of a real-valued
scoring function for which we formulate
 appropriate performance measures generalizing the AUC criterion.
We point out that ranking the best instances is an involved task
as it is a two-fold problem: (i) find the best instances and
(ii) provide a good ranking on these instances. The fact that
these two goals cannot be considered independently will be
highlighted in the paper. Despite this observation, we will first
formulate the issue of finding the best instances which is to be
understood as a toy problem for our purpose. This problem
corresponds to a {\em binary classification problem with a mass
constraint} (where the proportion $u_0$ of +1 labels predicted by
the classifiers is fixed) and it might present an interest {\em
per se}. The main complication here has to do with the necessity of
performing quantile estimation which affects the performance of
statistical procedures. Our proof technique was inspired by the former work
of Koul \cite{Kou02} in the context of $R$-estimation where similar statistics
arise.

\medskip

The rest of the paper is organized as follows. We first state the
problem of finding the best instances and study the performance of
empirical risk minimization in this setup (Section \ref{sec:cmc}).
We also explore the conditions on the distribution in order to
recover fast rates of convergence. In Section \ref{sec:crit} we
formulate performance measures for local ranking and provide
extensions of the AUC criterion. Eventually (Section
\ref{sec:erm}), we state some preliminary statistical results on
empirical risk minimization of these new criteria.

\section{Finding the best instances}\label{sec:cmc}

In the present section, we have a limited goal which is only to
determine the best instances without bothering of their order in
the list. By considering this subproblem, we will identify the
main technical issues involved in the sequel. It also permits to
introduce the main notations of the paper.

\medskip

Just as in standard binary classification, we consider the pair of
random variables $(X,Y)$ where $X$ is an observation vector in a
measurable space $\X$ and $Y$ is a binary label in $\{-1, +1\}$.
The distribution of $(X,Y)$ can be described by the pair $(\mu,
\eta)$ where $\mu$ is the marginal distribution of $X$ and $\eta$
is the a posteriori distribution defined by $\eta(x) = \Prob{Y=1
\mid X=x}$, $\forall x\in\X$. We define the {\em rate of best
instances} as the proportion of best instances to be considered
and denote it by $u_0\in(0,1)$. We denote by $Q(\eta, 1-u_0)$ the
$(1-u_0)$-quantile of the random variable $\eta(X)$. Then the {\em
set of best instances at rate $u_0$} is given by:
\[
C^*_{u_0} = \{x \in\X ~|~ \eta(x) \ge Q(\eta, 1-u_0)\}~.
\]

We mention two trivial properties of the set $C^*_{u_0}$ which
will be important in the sequel:
\begin{itemize}
\item {\sc Mass constraint:} we have $\mu\bigl( C^*_{u_0}
\bigr) = \Prob{X\in C^*_{u_0} }= u_0$,
\item {\sc Invariance property:} as a functional of
$\eta$, the set $C^*_{u_0}$ is invariant by strictly increasing
transforms of $\eta$.
\end{itemize}


\medskip

The problem of finding a proportion $u_0$ of the best instances
boils down to the estimation of the unknown set $C^*_{u_0}$ on the
basis of empirical data. Before turning to the statistical
analysis of the problem, we first relate it to binary
classification.

\subsection{A classification problem with a mass constraint}

A classifier is a measurable function $g ~:~ \X \to \{-1, +1\}$
and its performance is measured by the classification error $L(g)
= \Prob{Y\neq g(X)}$. Let $u_0\in(0,1)$ be fixed. Denote by
$g^*_{u_0} = 2 \I_{C^*_{u_0}}-1$ the classifier predicting +1 on
the set of best instances $C^*_{u_0}$ and -1 on its complement.
The next proposition shows that $g^*_{u_0}$ is an optimal element
for the problem of minimization of $L(g)$ over the family of
classifiers $g$ satisfying the {\em mass constraint}
$\Prob{g(X)=1}=u_0$.

\begin{proposition}
For any classifier $g~:~\X\to \{-1,+1\}$ such that $g(x) = 2
\I_C(x)-1$ for some subset $C$ of $\X$ and  $\mu(C) =
\Prob{g(X)=1}=u_0$, we have
\[
L^*_{u_0}\circeq L\bigl(g^*_{u_0}\bigr) \leq L(g)~.
\]
Furthermore, we have
\[
L^*_{u_0}=1-Q(\eta, 1-u_0)+(1-u_0)(2Q(\eta, 1-u_0)-1)-\EXP \left( \left| \eta(X)-Q(\eta,
1-u_0)\right|\right),
\]
and
\[
L(g)-L \bigl(g^*_{u_0}\bigr) =2\EXP \left( \left| \eta(X)-Q(\eta,
1-u_0)\right|\I_{C_{u_0}^*\Delta C}(X) \right),
\]
where $\Delta$ denotes the symmetric difference operation between
two subsets of $\X$.
\end{proposition}

\begin{proof}
For simplicity, we temporarily change the notation and set
$q=Q(\eta, 1-u_0)$. Then, for any classifier $g$ satisfying  the
the constraint $\Prob{g(X)=1}=u_0$, we have
\[
L(g) =  \EXP \left( (\eta(X)-q) \IND{g(X)=-1} + (q - \eta(X))
\IND{g(X)=+1} \right)  +  (1-u_0) q + (1-q) u_0 ~.
\]
The statements of the proposition immediately follow.\qed
\end{proof}

There are several progresses in the field of classification theory
where the aim is to introduce constraints in the classification
procedure or to adapt it to other problems. We relate our
formulation to other approaches in the following remarks.

\begin{remark}{\sc (Connection to hypothesis testing).}\label{rem:ht}
The implicit asymmetry in the problem due to the emphasis on the
best instances is reminiscent of the statistical theory of
hypothesis testing. We can formulate a test of simple hypothesis
by taking the {\em null assumption} to be $H_0 ~:~Y=+1$ and the
{\em alternative assumption} being $H_1 ~:~Y=-1$. We want to
decide which hypothesis is true given the observation $X$. Each
classifier $g$ provides a test statistic $g(X)$. The performance
of the test is then described by its type I error $\alpha(g) =
\Prob{g(X)=1 ~|~ Y=-1}$ and its power $\beta(g) = \Prob{g(X)=1 ~|~
Y=+1}$. We point out that if the classifier $g$ satisfies a mass
constraint, then we can relate the classification error with the
type I error of the test defined by $g$ through the relation:
\[
L(g) = 2(1-p) \alpha(g) +p-u_0
\]
where $p=\Prob{Y=1}$, and similarly, we have: $L(g) = 2p(1-
\beta(g)) - p-u_0$. Therefore, the optimal classifier minimizes
the type I error (maximizes the power) among all classifiers with
the same mass constraint. In some applications, it is more
relevant to fix a constraint on the probability of a false alarm
(type I error) and maximize the power. This question is explored
in a recent paper by Scott \cite{Sco06} (see also
\cite{ScoNow05}).
\end{remark}

\begin{remark}{\sc (Connection with regression level set estimation)}
We mention that the estimation of the level sets of the regression
function has been studied in the statistics literature
\cite{Cav97} (see also \cite{Tsy97}, \cite{WilNow06}) as well as
in the learning literature, for instance in the context of anomaly
detection (\cite{SHS05, ScoDav06,VerVer06}). In our framework of
classification with mass constraint, the threshold defining the
level set involves the quantile of the random variable $\eta(X)$.
\end{remark}

\begin{remark}{\sc (Connection with the minimum volume set approach)}
Although the point of view adopted in this paper is very
different, the problem described above may be formulated in the
framework of {\em minimum volume sets} learning as considered in
\cite{ScoNow06}. As a matter of fact, the set $C^*_{u_0}$ may be
viewed as the solution of the constrained optimization problem:
\[
\min_{C} \Prob{X\in C~|~ Y=-1}
\]
over the class of measurable sets $C$, subject to
\[
\Prob{X\in C} \ge u_0~.
\]
The main difference in our case comes from the fact that the
constraint on the volume set has to be estimated using the data
while in \cite{ScoNow06} it is computed from a known reference
measure. We believe that learning methods for minimum volume set
estimation may hopefully be extended to our setting. A natural way
to do it would consist in replacing conditional distribution of
$X$ given $Y=-1$ by its empirical counterpart. This is beyond the
scope of the present paper but will be the subject of future
investigation.
\end{remark}

\subsection{Empirical risk minimization}

We now investigate the estimation of the set $C_{u_0}^*$ of best instances at rate $u_0$ based on training data.
Suppose that we are given $n$ i.i.d. copies $(X_1,Y_1),\cdots, (X_n,Y_n)$ of the pair $(X,Y)$.
Since we have the ranking problem in mind, our methodology will consist in building the candidate sets from a
class $\S$ of real-valued scoring functions $s~:~\X\to \RR$. Indeed, we consider sets of the form
\[
C_s \circeq C_{s, u_0} = \{x \in\X ~|~ s(x) \ge Q(s, 1-u_0)\}~,
\]
where $s$ is an element of $\S$ and $Q(s, 1-u_0)$ is the $(1-u_0)$-quantile of the random variable $s(X)$. Note
that such sets satisfy the same properties of $C_{u_0}^*$ with respect to mass constraint and invariance to strictly
increasing transforms of $s$.

\medskip

From now on, we will take the simplified notation:
\[
L(s) \circeq L(s, u_0) \circeq L(C_{s}) = \Prob{Y \cdot (s(X) - Q(s, 1-u_0))<0}~.
\]

A scoring function minimizing the quantity
\[
L_n(s)=\frac{1}{n}\sum_{i=1}^n \I\{Y_i\cdot (s(X_i)-Q(s,1-u_0))<0\}.
\]
is expected to approximately minimize the true error $L(s)$, but the quantile depends on the unknown distribution of $X$.
In practice, one has to replace $Q(s,1-u_0)$ by its empirical counterpart $\hat{Q}(s,1-u_0)$ which corresponds to
the empirical quantile.
We will thus consider, instead of $L_n(s)$, the {\em truly empirical} error:
\[
\hat{L}_n(s)=\frac{1}{n}\sum_{i=1}^n \I\{Y_i\cdot (s(X_i)-\hat{Q}(s,1-u_0))<0\}.
\]
Note that $\hat{L}_n(s)$ is a complicated statistic since the empirical quantile involves all the instances
$X_1, \ldots, X_n$. We also mention that $\hat{L}_n(s)$ is a biased estimate of the classification error $L(s)$
of the classifier $g_s(x)=2\I\{s(x)\geq Q(s, 1-u_0)\}-1$.

\medskip

We introduce some more notations. Set, for all $t\in\RR$:
\begin{itemize}
\item $F_s(t) = \Prob{s(X)\le t}$
\item $G_s(t) = \Prob{s(X)\le t ~|~ Y=+1}$
\item $H_s(t) = \Prob{s(X)\le t ~|~ Y=-1}$
\end{itemize}
to be the cumulative distribution functions (cdf) of $s(X)$
(respectively, given $Y=1$, given $Y=-1$). We recall that the
definition of the quantiles of (the distribution of) a random
variable involves the notion of generalized inverse $F^{-1}$ of a
function $F$:
\[
F^{-1} (z) = \inf \{t\in\RR ~|~ F(t) \ge z \}~.
\]
Thus, we have, for all $v\in(0,1)$:
\[
Q(s,v) = F_s^{-1}(v) \quad \mbox{ and } \quad \hat{Q}(s,v)  = \hat{F}_s^{-1}(v)
\]
where $\hat{F}_s$ is the empirical cdf of $s(X)$: $\hat{F}_s(t) = \frac{1}{n}\sum_{i=1}^n \I\{s(X_i)\le t\}$, $\forall t\in\RR$.

\medskip

Without loss of generality, we will assume that all scoring functions in $\mathcal{S}$ take their values in
$(0,\lambda)$ for some $\lambda>0$. We now turn to study the performance of minimizers of $\hat{L}_n(s)$ over a class
$\S$ of scoring functions defined by
\[
\hat{s}_n=\argmin_{s\in \mathcal{S}}\hat{L}_n(s).
\]

Our first main result is an excess risk bound for the empirical risk minimizer $\hat{s}_n$ over a class $\S$
of uniformly bounded scoring functions. In the following theorem, we consider that the level sets of scoring
functions from the class $\S$ form a Vapnik-Chervonenkis (VC) class of sets.

\begin{theorem}\label{thm:firstorder}
We assume that
\begin{enumerate}
\item[(i)] the class $\S$ is symmetric (i.e. if $s\in\S$ then $\lambda-s\in\S$) and is a VC major
class of functions with VC dimension $V$.
\item[(ii)] the family $\K = \{~G_s, H_s ~:~s\in\S~\}$ of cdfs satisfies the following
property: any $K\in\K$ has left and right derivatives, denoted by
$K'_+$ and $K'_-$, and there exist strictly positive constants
$b$, $B$ such that $\forall (K,t)\in \K\times (0,\lambda)$,
\[
b \le \left|K'_+(t)\right| \le B \quad \mbox{ and } \quad b \le
\left|K'_-(t)\right| \le B ~.
\]
\end{enumerate}
For any $\delta>0$, we have, with probability larger than $1-\delta$,
\[
L(\hat{s}_n)-\inf_{s\in \S} L(s)\leq c_1 \sqrt{\frac{V}{n}} + c_2 \sqrt{\frac{ln(1/\delta)}{n}},
\]
for some positive constants $c_1, c_2$.
\end{theorem}

We now provide some insights on conditions (i) and (ii) of the theorem.

\begin{remark}{\sc (on the complexity assumption)}
On the terminology of major sets and major classes, we refer to Dudley \cite{Dud99}. In the proof, we need to
control empirical processes indexed by sets of the form $\{x ~:~ s(x) \ge t\}$ or $\{x ~:~ s(x) \le t\}$.
Condition (i) guarantees that these sets form a VC class of sets.
\end{remark}

\begin{figure}
\centering
\includegraphics[width=.5\textwidth]{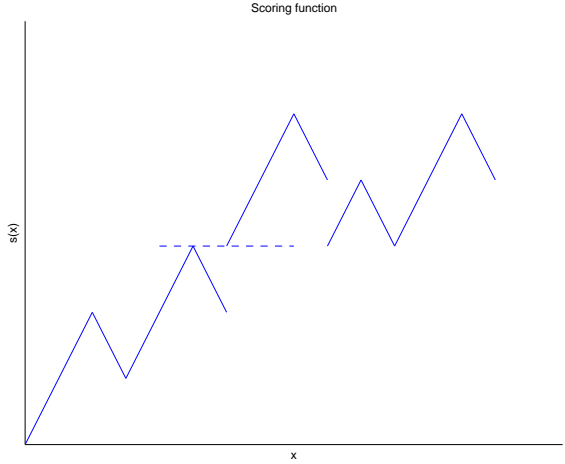}
\caption{Typical example of a scoring function.}
\label{fig:lbpscoring}
\end{figure}

\begin{remark}{\sc (on the choice of the class $\S$ of scoring functions)}
In order to grasp the meaning of condition (ii) of the theorem, we
consider the one-dimensional case with real-valued scoring
functions. Assume that the distribution of the random variable
$X_i$ has a bounded density $f$ with respect to Lebesgue measure.
Assume also that scoring functions $s$ are differentiable except,
possibly, at a finite number of points, and derivatives are
denoted by $s'$. Denote by $f_s$ the density of $s(X)$. Let $t\in
(0, \lambda)$ and denote by $x_1$, ..., $x_p$ the real roots of
the equation $s(x)=t$. We can express the density of $s(X)$ thanks
to the change-of-variable formula (see e.g. \cite{Pap65}):
\[
f_s(t) = \frac{f(x_1)}{s'(x_1)} + \ldots + \frac{f(x_p)}{s'(x_p)}~.
\]
This shows that the scoring functions should not present neither
flat nor steep parts. We can take for instance, the class $\S$ to
be the class of linear-by-parts functions with a finite number of
local extrema and with uniformly bounded left and right
derivatives: $\forall s\in\S$, $\forall x$, $m \le s'_+(x) \le M$
and $m \le s'_-(x) \le M$ for some strictly positive constants
$m$, and $M$ (see Figure \ref{fig:lbpscoring}). Note that any
subinterval of $[0, \lambda]$ has to be in the range of scoring
functions $s$ (if not, some elements of $\K$ will present a {\em
plateau}). In fact, the proof requires such a behavior only in the
vicinity of the points corresponding to the quantiles $Q(s,
1-u_0)$ for all $s\in\S$.
\end{remark}

\begin{proof}
Set $v_{0}=1-u_{0}$. By a standard argument (see e.g. \cite{DGL96}), we have:
\begin{align*}
L(\hat{s}_n)-\inf_{s\in \S}L(s) & \le  2\sup_{s\in \S}\left|\hat{L}_n(s)-L(s)\right|\\
& \le 2 \sup_{s\in \S}\left|\hat{L}_n(s)-L_n(s)\right| + 2 \sup_{s\in \S}\left|L_n(s)-L(s)\right|~.
\end{align*}

Note that the second term in the bound is an empirical process whose behavior is well-known. In our case, assumption (i)
implies that the class of sets $\{x \, : s(x)\ge Q(s, v_0) \}$ indexed by scoring functions $s$ has a VC dimension smaller than $V$.
Hence, we have by a concentration argument combined with a VC bound for the expectation of the supremum
(see, e.g. \cite{Lug02}), for any $\delta>0$, with probability larger than $1-\delta$,
\[
\sup_{s\in \mathcal{S}}\left|L_n(s)-L(s)\right| \le c\sqrt{\frac{V}{n}} + c'\sqrt{\frac{ln(1/\delta)}{n}}
\]
for universal constants $c, c'$.

\medskip

We now show how to handle the first term. Following the work of
Koul \cite{Kou02}, we set the following notations:
\begin{align*}
M(s,v) = & ~\Prob{ Y \cdot \bigl(s(X)-Q(s,v)\bigr)<0} \\
U_n(s,v) =& ~\frac{1}{n}\sum_{i=1}^n \I\{Y_i \cdot
\bigl(s(X_i)-Q(s,v)\bigr)<0\} -M(s,v) ~.
\end{align*}
and note that $U_n(s,v)$ is centered.

\medskip

We then have the following decomposition, for any $s\in \S$ and $v_0\in (0,1)$:
\[
\left|\hat{L}_n(s)-L_n(s)\right| \le \left|U_n(s, F_s\circ\hat{F}_s^{-1}(v_0))-U_n(s,v_0)\right|
+ \left|M(s, F_s\circ\hat{F}_s^{-1}(v_0))-M(s,v_0)\right|~.
\]

Note that $M(s, F_s\circ\hat{F}_s^{-1}(v_0)) = \Prob{ Y \cdot \bigl(s(X)-\hat{Q}(s,v)\bigr)<0 \mid D_n}$
where $D_n$ denotes the sample $(X_1,Y_1),\cdots, (X_n,Y_n)$.

\medskip

Recall the notation $p=\Prob{Y=1}$. Since $M(s, v)=(1-p)(1-H_s\circ F_s^{-1}(v))+p G_s\circ F_s^{-1}(v)$ and $F_s=p G_s + (1-p) H_s$,
the mapping $v\mapsto M(s,v)$ is Lipschitz by assumption (ii).
Thus, there exists a constant $\kappa<\infty$, depending only on $p$, $b$ and $B$, such that:
\[
\left|M(s,F_s\circ\hat{F}_s^{-1}(v_0))-M(s,v_0)\right|\le \kappa \left|F_s\circ\hat{F}_s^{-1}(v_0)-v_0)\right|.
\]

Moreover, we have, for any $s\in\S$:
\begin{align*}
\left|F_s\circ\hat{F}_s^{-1}(v_0)-v_0\right|  \le & \left|F_s\circ\hat{F}_s^{-1}(v_0)-\hat{F}_s\circ\hat{F}_s^{-1}(v_0))\right|
+\left|\hat{F}_s\circ\hat{F}_s^{-1}(v_0)-v_0)\right|\\
\le & \sup_{t\in (0,\lambda)} \left| F_s(t)-\hat{F}_s(t) \right|+ \frac{1}{n}~.
\end{align*}

Here again, we can use assumption (i) and a classical VC bound from \cite{Lug02} in order to
control the empirical process, with probability larger than $1-\delta$:
\[
\sup_{(s,t) \in \S \times (0,\lambda)} \left| F_s(t)-\hat{F}_s(t) \right| \le c\sqrt{\frac{V}{n}} + c' \sqrt{\frac{ln(1/\delta)}{n}}
\]
for some constants $c, c'$.

\medskip

It remains to control the term involving the process $U_n$:
\begin{align*}
\left|U_n(s, F_s\circ\hat{F}_s^{-1}(v_0))-U_n(s,v_0)\right| &\le \sup_{v \in (0,1)}\left|U_n(s, v)-U_n(s,v_0)\right| \\
 & \le 2 \sup_{v \in (0,1)}\left|U_n(s, v)\right|
\end{align*}

Using that the class of sets of the form $\{x \, : s(x)\ge Q(s, v) \}$ for  $v\in (0,1)$ is included
in the class of sets of the form $\{x \, : s(x)\ge t \}$ where $t\in (0,\lambda)$, we then have
\[
\sup_{v \in (0,1)}\left|U_n(s, v)\right| \le \sup_{t\in (0,\lambda)} \left|\frac{1}{n}\sum_{i=1}^n \I\{Y_i \cdot \bigl(s(X_i)-t\bigr)<0\} - \Prob{ Y \cdot \bigl(s(X)-t\bigr)<0} \right|
 ~,
\]
which leads again to an empirical process indexed by a VC class of sets and can be bounded as before.\qed
\end{proof}

\subsection{Fast rates of convergence}

We now propose to examine conditions leading to fast rates of
convergence (faster than $n^{-1/2}$). It has been noticed (see
\cite{MamTsy99}, \cite{Tsy04}, \cite{MasNed06}) that it is
possible to derive such rates of convergence in the classification
setup under additional assumptions on the distribution. We propose
here to adapt these assumptions for the problem of classification
with mass constraint.

\medskip

Our concern here is to formulate the type of conditions which render the problem easier from a statistical perspective.
For this reason and to avoid technical issues, we will consider a quite restrictive setup where it is assumed that:
\begin{enumerate}
\item the class $\S$ of scoring functions is a finite class with $N$ elements,
\item an optimal scoring rule $s^*$ is contained in $\S$.
\end{enumerate}

We have found that the following additional conditions on the
distribution and the class $\S$ allow to derive fast rates of
convergence for the excess risk in our problem.

\begin{itemize}
\item[(iii)] There exist constants $\alpha \in (0,1)$ and $B>0$ such that, for all $t\ge 0$,
\[
\Prob{\left|  \eta(X)-Q(\eta, 1-u_0) \right| \le t} \le B \, t^{\frac{\alpha}{1-\alpha}}~.
\]
\item[(iv)] the family $\K = \{~G_s, H_s ~:~s\in\S~\}$ of cdfs satisfies the following
property: for any $s\in \S$, $G_s$ and $H_s$ are twice differentiable at $Q(s,1-u_0)=F_s^{-1}(1-u_0)$.
\end{itemize}

We point out that conditions (ii) and (iii) are not completely
independent. Indeed, if $(G_{\eta}, H_{\eta})$ belongs to the
class $\K$ fulfilling condition (ii), then
$F_{\eta}=pG_{\eta}+(1-p)H_{\eta}$ is Lipschitz and condition
(iii) is satisfied with $\alpha= 1/2$. Note that condition (iii)
simply extends the standard low noise assumption introduced by
Tsybakov \cite{Tsy04} (see also \cite{BBL05} for an account on
this) where the level 1/2 is replaced by the $(1-u_0)$-quantile of
$\eta(X)$. Indeed, we have, under condition (iii), the variance
control, for any $s\in\S$:
\[
\Var (\I\{Y\neq 2\I_{C_{s}}(X)-1\}-\I\{Y\neq
2\I_{C^*_{u_0}}(X)-1\}) \le c~(L(s)-L^*_{u_0})^{\alpha}~,
\]
or, equivalently,
\[
\EXP \bigl(\I_{C_{s}\Delta C^*_{u_0}}(X)\bigr) \le
c~(L(s)-L^*_{u_0})^{\alpha}~.
\]

Now, if we denote
\[
s_n=\argmin_{s\in \S} L_n(s)~,
\]
we have, by a standard argument based on Bernstein's inequality
(see Section 5.2 in \cite{BBL05}), with probability $1-\delta$,
\[
L(s_n)-L^*_{u_0}\leq c \, \left(\frac{log(N/\delta)}{n}\right)^{\frac{1}{2-\alpha}}~.
\]
for some positive constant $c$.

\medskip

The novel part of the analysis below lies in the control of the
bias induced by plugging the empirical quantile $\hat{Q}(s,
1-u_0)$ in the risk functional. The next theorem shows that faster
rates of convergence can be obtained under the previous
assumptions with $\alpha=1/2$.

\begin{theorem} \label{thm:fast}
We assume that the class $\S$ of scoring functions is a finite
class with $N$ elements, and that it contains an optimal scoring
rule $s^*$. Moreover, we assume that conditions (i)-(iv) are
satisfied. Then, for any $\delta>0$, we have, with probability
$1-\delta$:
\[
L(\hat{s}_n)-L^*_{u_0}  \le \displaystyle c \, \left(\frac{log(N/\delta)}{n}\right)^{\frac{2}{3}},
\]
for some constant $c$.
\end{theorem}

\begin{remark}{\sc (on the rate $n^{-2/3}$)}
The previous results highlights the fact that rates faster than
the one obtained in Theorem \ref{thm:firstorder} can be obtained
in this setup with additional regularity assumptions. However, it
is noteworthy that the standard low noise assumption (iii) is
already contained in assumption (ii) which is required in proving
the typical $n^{-1/2}$ rate. The consequence of this observation
is there is no hope of getting rates up to $n^{-1}$ unless
assumption (ii) is weakened.
\end{remark}


The proof of the previous theorem is based on two arguments: the
structure of {\em linear signed rank statistics} and the variance
control assumption. The situation is similar to the one we
encoutered in \cite{CLV06} where we were dealing with {\sc
U}-statistics and we had to invoke Hoeffding's decomposition in
order to grasp the behavior of the underlying {\sc U}-processes.
Here we require a similar argument to describe the structure of
the empirical risk functional $\hat{L}_n(s)$ under study. This
statistic can be interpreted as a linear signed rank statistic and
the key decomposition has been used in the context of
nonparametric hypotheses testing and {\sc R}-estimation. We mainly
refer to H\'ajek and Sidak \cite{HajSid67}, Dupac and H\'ajek
\cite{DupHaj69}, Koul \cite{Kou70}, Koul and Staudte
\cite{KouSta72} for an account on rank statistics.

\medskip

We briefly go through the main ideas, but first we need to
introduce some notations. Set:
\begin{align*}
\forall v\in[0, 1]~, \quad K(s,v)&=\EXP \left( Y\I \{s(X)\le Q(s,v)\} \right) = p G_s(Q(s,v))-(1-p) H_s(Q(s,v))\\
\hat{K}_n(s,v)&= \frac{1}{n} \sum_{i=1}^n Y_i\I \{s(X_i)\le
\hat{Q}(s,v)\}~.
\end{align*}

Then we can write:
\begin{align*}
L(s)&= 1-p + K(s,1-u_0)\\
\hat{L}_n(s)&=  \frac{n_-}{n}+\hat{K}_n(s,1-u_0)~,
\end{align*}
where $n_-=\sum_{i=1}^n \I \{Y_i=-1\}$.

\medskip

We note that the statistic $\hat{L}_n(s)$ is related to linear
signed rank statistics.

\begin{definition}{[Linear signed rank statistic].}
Consider $Z_1, \ldots, Z_n$ an i.i.d. sample with distribution $F$
and a real-valued score generating function $\Phi$. Denote by
$R^+_i = \mbox{rank}(\vert Z_i \vert)$ the rank of $\vert Z_i
\vert$ in the sample $\vert Z_1 \vert, \ldots, \vert Z_n \vert$.
Then the statistic
\[
\sum_{i=1}^{n} \Phi\left(\frac{R^+_i}{n+1}\right)~\sgn(Z_i)
\]
is a linear signed rank statistic.
\end{definition}

\begin{proposition}
For fixed $s$ and $v$, the statistic $\hat{K}_n(s, v)$ is a linear
signed rank statistic.
\end{proposition}

\begin{proof}
Take $Z_i=Y_is(X_i)$. The random variables $Z_i$ have their
absolute value distributed according to $F_s$ and have the same
sign as $Y_i$. It is easy to see that the statistic $\hat{K}_n(s,
v)$ is a linear signed rank statistic with score generating
function $\Phi(x) = \I_{[x\le v]}$. \qed
\end{proof}

A decomposition of Hoeffding's type for such statistics can be
formulated. Set first:
\[
Z_n(s, v) = \frac{1}{n} \sum_{i=1}^n \left( Y_i - K'(s,v)
\right)\I\{s(X_i)\le Q(s,v)\}-K(s,v)+v K'(s,v)~,
\]
where  $K'(s,v)$ denotes the derivative of the function $v \mapsto
K(s, v)$. Note that $Z_n(s, v)$ is a centered random variable with
variance:
\[
\sigma^2(s,v)=v-K(s,v)^2 + v(1-v) K'^2(s,v) -2 (1-v)
K'(s,v)K(s,v)~.
\]

\medskip

The next result is due to Koul \cite{Kou70} and we provide an
alternate proof in the Appendix.

\begin{proposition}\label{prop:hoeff}
We have, for all $s\in \S$ and $v\in[0, 1]$:
\[
\hat{K}_n(s,v)=K(s,v)+Z_n(s, v)+\Lambda_n(s)~.
\]
with
\[
\Lambda_n(s) = O_{\PROB}(n^{-1}) ~ \mbox{ as } ~ n\rightarrow
\infty~.
\]
\end{proposition}

This asymptotic expansion highlights the structure of the
statistic $L_n(s)$ for fixed $s$. The leading term $Z_n(s, 1-u_0)$
is an empirical average of i.i.d. random variables and it provides
the asymptotic variance of $L_n(s)$. It is worth noticing that
$Z_n(s, 1-u_0)$ is not reduced to the same empirical functional
with the true, instead of the empirical, quantile but it also
involves a derivative term. Since the remainder term
$\Lambda_n(s)$ is of the order $n^{-1}$, it will not affect the
final rate of convergence under low noise conditions. Therefore,
the variance control assumption concerns the variance of the
function involved in the empirical average $Z_n(s, 1-u_0)$.

\medskip

We denote by:
\[
h_s(X_i, Y_i) = \left( Y_i - K'(s,v) \right)\I\{s(X_i)\le Q(s,v)\}
- K(s,v)+ v K'(s,v)~,
\]
and we then have
\[
Z_n(s, v) -Z_n(s^*, v)  = \frac{1}{n} \sum_{i=1}^n \bigl(h_s(X_i,
Y_i) - h_{s^*}(X_i, Y_i) \bigr)  ~.
\]

\begin{proposition}\label{prop:varcont}
Fix $v\in [0, 1]$. Assume that condition (iii) holds. Then, we
have, for all $s\in \S$:
\[
\Var \bigl(h_s(X_i, Y_i) - h_{s^*}(X_i, Y_i) \bigr) \le c
\bigl(L(s) - L(s^*) \bigr)^{\alpha} ~,
\]
for some constant $c$.
\end{proposition}

\begin{proof}
We first write that:
\[
h_s(X_i, Y_i) - h_{s^*}(X_i, Y_i) = I + II + III + IV + V
\]
where
\[
\begin{array}{ccl}
I & = & Y_i \ \bigl(\I\{s(X_i)\le Q(s,v)\} - \I\{s^*(X_i)\le Q(s^*,v)\} \bigr) \\
&&\\
II & = & (K'(s^*,v) - K'(s,v))~\I\{s^*(X_i)\le Q(s^*,v)\} \\
&&\\
III & = & K'(s,v) ~ \bigl(\I\{s^*(X_i)\le Q(s^*,v)\} - \I\{s(X_i)\le Q(s,v)\} \bigr)  \\
&&\\
IV & = & K(s^*, v) - K(s, v)\\
&&\\
V & = & v ~(K'(s, v) - K'(s^*,v))~.
\end{array}
\]

\smallskip

By Cauchy-Schwarz inequality, we only need to show that the
expected value of the square of these quantities is smaller than
$(L(s)-L^*)^{\alpha}$ up to some multiplicative constant.

\medskip

Note that, by definition of $K$, we have:
\[
\begin{array}{ccl}
II & = &  (L'(s^*, v) - L'(s, v))~ \I\{s^*(X_i)\le Q(s^*,v)\} \\
&&\\
IV & = & L(s^*) - L(s)\\
&&\\
V & = &  v ~ (L'(s, v) - L'(s^*, v))
\end{array}
\]
where $L'(s,v)$ denotes the derivative of the function $v \mapsto
L(s, v)$. It is clear that, for any $s$, we have $L(s, v) = L(s^*,
v)$ implies that $L'(s, v) = L'(s^*, v)$ otherwise $s^*$ would not
be an optimal scoring function at some level $v'$ in the vicinity
of $v$. Therefore, since $\S$ is finite, there exists a constant
$c$ such that
\[
(L'(s, v) - L'(s^*, v))^2 \le c (L(s)-L^*)^{\alpha}
\]
and then $\EXP\bigl( II^2 \bigr)$ and $\EXP\bigl( V^2 \bigr)^2$
are bounded accordingly.

\medskip

Moreover, we have:
\begin{align*}
\EXP\bigl( I^2 \bigr) & \le \EXP \bigl(\I_{C_s
\Delta C_{s^*}}(X) \bigr) \\
& \\
& \le c(L(s) - L(s^*))^{\alpha}
\end{align*}
for some positive constant $c$, by assumption (iii).

\medskip

Eventually, by assumption (ii), we have that $K'(s,v)$ is
uniformly bounded and thus, the term $\EXP\bigl( III^2 \bigr)$ can
be handled similarly.
 \qed
\end{proof}

\bigskip

\noindent {\bf Proof of Theorem \ref{thm:fast}.} The proof is the
same as the one of Theorem 5 from \cite{CLV06} which uses a result
by Massart \cite{Mas06}. \qed

\section{Performance measures for local ranking}\label{sec:crit}

Our main interest here is to develop a setup describing the
problem of not only finding but also ranking the best instances.
As far as we know, this problem has not been considered from a
statistical perspective  until now. In the sequel, we build on the
results from Section \ref{sec:cmc} and also on our previous work
on the (global) ranking problem \cite{CLV06} in order to capture
some of the features of the local ranking problem. The present
section is devoted to the construction of performance measures
reflecting the quality of ranking rules on a restricted set of
instances.

\subsection{ROC curves and optimality in the local ranking problem}

We consider the same statistical model as before with $(X, Y)$
being a pair of random variables over $\X \times \{-1, +1\}$ and
we examine ranking rules resulting from real-valued scoring
functions $s~:~\X\to (0, \lambda)$. The reference tool for
assessing the performance of a scoring function $s$ in separating
the two populations (positive vs. negative labels) is the
Receiving Operator Characteristic known as the ROC curve
(\cite{vTr68}, \cite{Ega75}). If we take the notations
$\bar{G}_s(z) = \Prob{s(X)>z ~\mid~Y=1}$ (true positive rate) and
$\bar{H}_s(z) = \Prob{s(X)>z ~\mid~Y=-1}$ (false positive rate),
we can define the ROC curve, for any scoring function $s$, as the
plot of the function:
\[
z \mapsto \left(\bar{H}_s(z), \bar{G}_s(z) \right)
\]
for thresholds $z\in(0, \lambda)$, or equivalently as the plot of the function:
\[
t \mapsto \bar{G}_s \circ H^{-1}_s(1-t)
\]
for $t\in (0,1)$. The optimal scoring function is the one whose
ROC curve dominates all the others for all $z\in (0, \lambda)$ (or
$t\in (0,1)$) and such a function actually exists. Indeed, by
recalling the hypothesis testing framework in the classification
model (see Remark \ref{rem:ht}) and using Neyman-Pearson's Lemma,
it is easy to check that ROC curve of the function $\eta(x) =
\Prob{Y=1~\mid~X=x}$ dominates the ROC curve of any other scoring
function. We point out that the ROC curve of a scoring function
$s$ is invariant by strictly increasing transformations of $s$.

\medskip

In our approach, for a given scoring function $s$, we focus on thresholds $z$ corresponding to the cut-off separating a proportion $u\in (0,1)$ of
top scored instances according to $s$ from the rest. Recall from Section \ref{sec:cmc} that the best instances according to $s$ are the elements of the set
$C_{s, u} =  \{x \in\X ~|~ s(x) \ge Q(s, 1-u)\}$ where $Q(s, 1-u)$ is the $(1-u)$-quantile of $s(X)$. We set the following notations:
\begin{align*}
\alpha (s, u) & = \Prob{s(X) \ge Q(s,1-u) \mid Y=-1}\\
\beta (s, u) & = \Prob{s(X) \ge Q(s,1-u) \mid Y=+1}~.
\end{align*}

We propose to re-parameterize the ROC curve with the proportion $u\in (0,1)$ and then describe it as the plot of the function:
\[
u \mapsto (\alpha (s, u), \beta (s, u))~,
\]
for each scoring function $s$. When focusing on the best instances at rate $u_0$, we only consider
the part of the ROC curve for values $u\in (0, u_0)$.

\medskip

However attractive is the ROC curve as a graphical tool, it is not
a practical one for developing learning procedures achieving
straightforward optimization. The most natural approach is to
consider risk functionals built after the ROC curve such as the
Area Under an ROC Curve (known as the AUC or AROC, see
\cite{HaMcNe82}). Our goals in this section are:
\begin{enumerate}
\item to extend the AUC criterion in order to focus on restricted parts of the ROC curve,
\item to describe the optimal elements with respect to this extended criterion.
\end{enumerate}

\begin{figure}
\centering
\includegraphics[width=6.3cm]{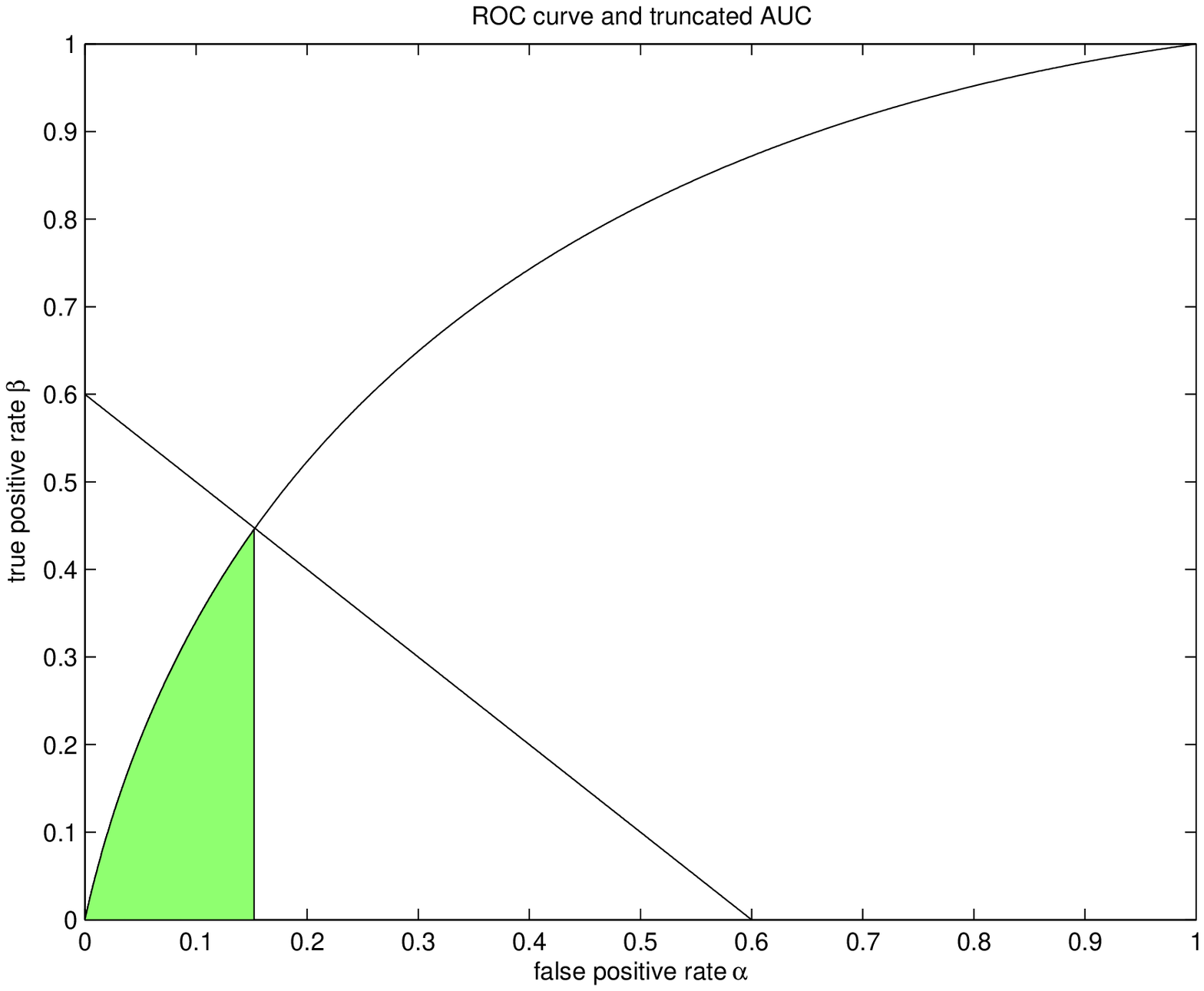}\hfill
\includegraphics[width=6.3cm]{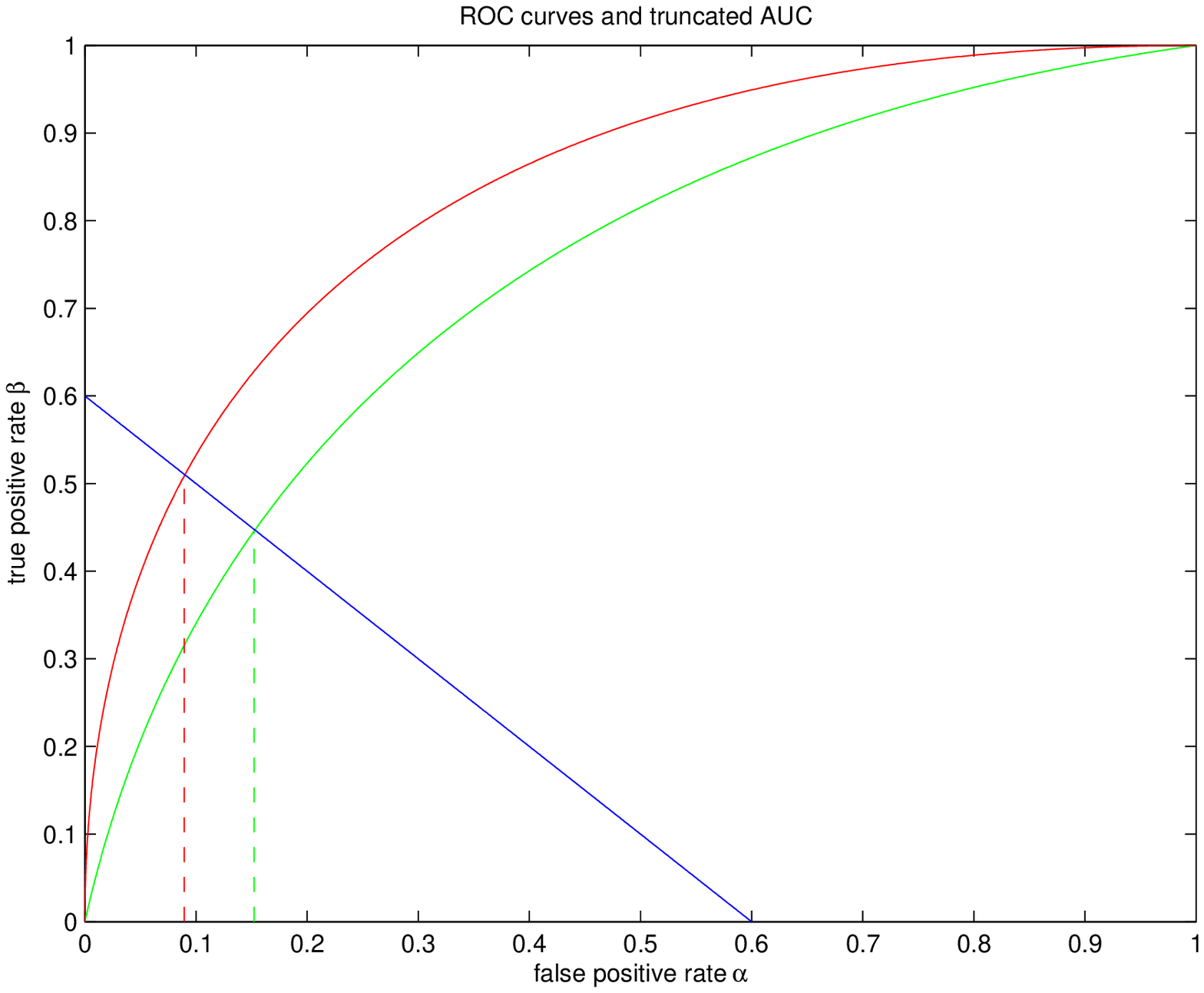}
\caption{ROC curves, line $D(u_0, p)$ and truncated AUC at rate
$u_0$ of best instances. } \label{fig:truncauc}
\end{figure}

We point out the fact that extending the AUC is not trivial.
Indeed, we notice that $\alpha (s, u)$ and  $\beta (s, u)$ are
related by a linear relation, for fixed $u$ and $p$, when $s$
varies:
\[
u = p\beta (s, u) + (1-p) \alpha (s, u)
\]
where $p =\Prob{Y=1}$. We denote the line plot of this relation by
$D(u, p)$. Hence, the part of the ROC curve of a scoring function
$s$ corresponding to the best instances at rate $u_0$ is the part
going from the origin $(0,0)$ to the intersection between the line
$D(u_0, p)$ and the ROC curve (shaded area in the left display of
Figure \ref{fig:truncauc}). It follows that, the closer to $\eta$
the scoring function $s$ is, the higher the ROC curve is, but at
the same time the integration domain shrinks (right display of
Figure \ref{fig:truncauc}).

\medskip

Our guideline in defining risk criteria for the problem of ranking
the best instances is the form of the optimal elements. We expect
the optimal scoring functions at the rate $u_0$ to belong to the
equivalence class (functions defined up to the composition with a
nondecreasing transformation) defined by scoring functions $s^*$
such that:
\[
s^*(x) = \left\{
\begin{array}{lll}
 & \eta(x) & \qquad \mbox{ if } x\in C^*_{u_0} \\
\\
< & \displaystyle  \inf_{C^*_{u_0}} \eta & \qquad \mbox{ if }
x\notin C^*_{u_0}~.
\end{array}
\right.
\]

Such scoring functions fulfill the two properties of finding the
best instances (indeed $C_{s^*, u_0} = C^*_{u_0}$) and ranking
them as well as the regression function. We will denote by $\S^*$
the set of optimal scoring functions for the problem of ranking
the best instances at the rate $u_0$.

\medskip

As a preliminary result, and before proposing an adequate criterion, we formulate a simple lemma.

\begin{lemma}\label{lem:opt}
For any scoring function $s$, we have for all $u\in (0,1)$,
\begin{align*}
\beta (s, u) & \le \beta (\eta, u) \\
\alpha (s, u) & \ge \alpha (\eta, u)~.
\end{align*}
Moreover, we have equality only for those $s$ such that $C_{s, u_0} = C^*_{u_0}$.
\end{lemma}

\begin{proof}
We show the first inequality. By definition, we have:
\[
\beta (s, u) = 1-H_{s}(Q(s, 1-u))~.
\]
Observe that, for any scoring function $s$,
\begin{align*}
p(1-H_{s}(Q(s, 1-u)) & = \Prob{Y=1, s(X)>Q(s, 1-u)} \\
& = \EXP \left(\eta(X)\I\{X\in C_{s,u}\}\right)~.
\end{align*}

We thus have
\begin{multline*}
p \left( H_{s}(Q(s, 1-u) - H_{\eta}(Q(\eta, 1-u)\right) \\
\\
\qquad
\begin{array}{l}
= \EXP \left(\eta(X)(\I\{X\in C^*_{u}\}-\I\{X\in C_{s,u}\}) \right)\\
\\
= \EXP \left( \eta(X)\I\{X\notin C^*_{u}\}\left(\I\{X\in C^*_{u}\}-\I\{X\in C_{s,u}\}\right)\right) \\
\\
\qquad  +\EXP\left(\eta(X)\I\{X\in C^*_{u}\}(\I\{X\in C^*_{u}\}-\I\{X\in C_{s,u}\})\right)\\
\\
\geq  - \EXP\left(Q(\eta, 1-u)\I\{X\notin C^*_{u}\}~\I\{X\in C_{s,u}\}\right)
+\EXP\left(Q(\eta, 1-u)\I\{X\in C^*_{u}\}(1-\I\{X\in C_{s,u}\})\right)\\
\\
= Q(\eta, 1-u)(1-u-1+u)=0~.
\end{array}
\end{multline*}

The second inequality simply follows from the identity below:
\[
1-u=p H_{s}(Q(s, 1-u))+(1-p)G_s(Q(s, 1-u))~. \mbox{ \qed}
\]
\end{proof}

\medskip

In view of this result, a wide collection of criteria with the set
$\S^*$ as the set of optimal elements could naturally be
considered, depending on how one wants to weight the two types of
error $1-\beta(s,u)=$ (type II error in the hypothesis testing
framework) and $\alpha(s,u)$ (type I error) according to the rate
$u\in[0,u_{0}]$. However, not all the criteria obtained in this
manner can be interpreted as generalizations of the AUC criterion
for $u_{0}=1$.

\subsection{Generalization of the {\sc AUC} criterion}

In \cite{CLV06}, we have considered the ranking error of a scoring function $s$
as defined by:
\[
R(s) = \PROB\{(Y-Y')(s(X)-s(X'))<0\}~,
\]
where $(X', Y')$ is  an i.i.d. copy of the random pair $(X, Y)$.

\medskip

Interestingly, it can be proved that minimizing the ranking error
$R(s)$ is equivalent to maximizing the well-known {\sc AUC}
criterion. This is trivial once we write down the probabilistic
interpretation of the {\sc AUC}:
\[
\mbox{{\sc AUC}}(s) =  \Prob{ s(X) > s(X') \mid Y=1,
\, Y'=-1} = 1-\frac{1}{2p(1-p)} R(s)~.
\]

We now propose a local version of the ranking error on a measurable set $C\subset
\X$:
\[
R(s, C)=\Prob{(s(X)-s(X'))(Y-Y')>0, ~ (X,X')\in C^2}~,
\]
and the local analogue of the AUC criterion:
\[
\mbox{\sc LocAUC}(s, u)=\Prob{s(X)>s(X'), ~s(X)\ge
Q(s,1-u)~\mid~ Y=1, Y'=-1}~.
\]

This criterion obviously boils down to the standard criterion for
$u=1$. However, in the case where $u<1$, we will see that there is
no equivalence between maximizing the {\sc LocAUC} criterion and
minimizing the local ranking error $s \mapsto R(s, u) \circeq R(s,
C_{s, u})$. Indeed, the local ranking error is not a relevant
performance measure for finding the best instances. Minimizing it
would solve the problem of finding the instances that are the
easiest to rank.

\medskip

The following theorem states that  scoring functions $s^*$ in the
set $\S^*$ maximize the  criterion {\sc LocAUC} and that the
latter may be decomposed as a sum of a 'power' term and (the
opposite of) a local ranking error term.

\begin{theorem}\label{thm:optauc}
Let $u_0\in (0,1)$. We have, for any scoring function $s$:
\[
\forall s^*\in\S^*, \quad  \mbox{\sc LocAUC} (s, u_0) \le \mbox{\sc
LocAUC}(s^*, u_0)~.
\]
Moreover, the following relation holds:
\[
\forall s, \quad  \mbox{\sc LocAUC}(s, u_0)
  = \beta(s,u_0)-\frac{1}{2p(1-p)}R(s, u_0)~,
\]
where $R(s, u_0) = R(s, C_{s, u_0})$.
\end{theorem}

\begin{proof}
Set $v_0=1-u_0$. Observe first that:
\begin{align*}
\mbox{\sc LocAUC}(s, u_{0})& =  \EXP \left(H_{s}(s(X))~\I\{s(X)\ge Q(s,v_{0})\} ~\mid Y=1~\right)\\
&\\
& = \int_{Q(s,v_{0})}^{+\infty} H_{s}(z)~G_{s}(dz)~.
\end{align*}

We use that $p G_s = F_{s}-(1-p)H_{s}$ and we obtain:
\begin{align*}
p\mbox{\sc LocAUC}(s, u_{0}) & = \int_{Q(s,v_{0})}^{+\infty} H_{s}(z)~F_{s}(dz) - (1-p)\int_{Q(s,v_{0})}^{+\infty} H_{s}(z)~H_{s}(dz)\\
&\\
&= \int_{v_{0}}^1 (1 - \alpha(s, v))~dv - \frac{1-p}{2}
\left(1-(1-\alpha(s, v_0))^2\right).
\end{align*}

This formula, combined with Lemma \ref{lem:opt}, establishes the first part of Theorem \ref{thm:optauc}.

\medskip

Besides, integrating by parts and making a change of variables, we get:
\begin{align*}
\int_{Q(s,v_{0})}^{+\infty} H_{s}(z)~G_{s}(dz) & =
1-(1-\alpha(s,u_{0}))(1-\beta(s,u_{0}))
- \int_{0}^{\alpha(s,u_{0})}(1-\beta(s,\alpha))~d\alpha \\
& \\
&= \int_{0}^{\alpha(s,u_{0})}\beta(s,\alpha) \,
d\alpha+\beta(s,u_{0})(1-\alpha(s,u_{0}))~.
\end{align*}

On the other hand, one has
\begin{align*}
\alpha(s,u_{0})\beta(s,u_{0})
& = \displaystyle \frac{1}{p(1-p)}\Prob{s(X) \wedge s(X')>Q(s, v_{0}) , ~ Y'=1 , ~ Y=-1} \\
\\
&= \displaystyle \Prob{s(X')>s(X) , ~s(X)\wedge s(X')>Q(s, v_{0})~\mid~ Y'=1 , ~ Y=-1} \\
\\
&\qquad \displaystyle + \frac{1}{p(1-p)}\Prob{s(X')<s(X) , ~ (X,X')\in C_{s,u_{0}}^2 , ~ Y'=1 , ~ Y=-1} \\
\\
&= \displaystyle\int_{0}^{\alpha(s,u_{0})}\beta(s,\alpha)~d\alpha
+ \frac{1}{2p(1-p)}R(s, u_{0})~.
\end{align*}
Plugging this in the previous formula leads to the second
statement of the theorem. \qed
\end{proof}

\begin{remark} {\sc (Truncating the AUC)}
In the theorem, we obviously recover the relation between the
standard {\sc AUC} criterion and the (global) ranking error when
$u_0=1$. Besides, by checking the proof, one may relate the
generalized {\sc AUC} criterion to the truncated {\sc AUC}. As a
matter of fact, we have:
\[
\forall s~, \quad \mbox{\sc LocAUC}(s,
u_{0})=\int_{0}^{\alpha(s,u_0)}\beta(s,
\alpha)~d\alpha+\beta(s,u_0)-\alpha(s,u_0)\beta(s,u_0).
\]
The values $\alpha(s,u_0)$ and $\beta(s,u_0)$ are the coordinates
of the intersecting point between the ROC curve of the scoring
function $s$ and the line $D(u_0, p)$. Thus, the integral term
represents the area of the surface delimited by the ROC curve, the
horizontal $x$-axis and the line $x=\alpha(s,u_0)$ (see Figure
\ref{fig:truncauc}). The theorem reveals that evaluating the local
performance of a scoring statistic $s(X)$ by the truncated AUC as
proposed in \cite{DodPep03} is highly arguable since the maximizer
of the functional $s \mapsto \int_{0}^{\alpha(s,u_0)}\beta(s,
\alpha)~d\alpha$ is usually not in $\S^*$.
\end{remark}

\subsection{Generalized Wilcoxon statistic}

We now propose a different extension of the plain {\sc AUC} criterion.
Consider $(X_1, Y_1)$, $\ldots$, $(X_n, Y_n)$, $n$ i.i.d. copies of the random pair $(X, Y)$.
The intuition relies on a well-known relationship
between Mann-Whitney and Wilcoxon statistics. Indeed, a natural empirical estimate of the {\sc AUC} is the {\em rate of concording pairs}:
\[
\widehat{\mbox{\sc AUC}}(s)=\frac{1}{n_+n_-}\sum_{1\le i, j \le n}\I\{Y_{i}=-1,Y_{j}=1, s(X_i)<s(X_j)\}~,
\]
with $n_+=n-n_-=\sum_{i=1}^n \I\{Y_i=+1\}$. On the other hand, we recall that the Wilcoxon statistic $T_{n}(s)$ is the two-sample linear
rank statistic associated to the {\em score generating function} $\Phi(v)=v$, $\forall v\in (0,1)$, obtained by summing the ranks corresponding to positive labels:
\[
T_{n}(s)=\sum_{i=1}^{n} \I\{Y_{i}=1\}~\frac{\mbox{rank}(s(X_{i}))}{n+1},
\]
where $\mbox{rank}(s(X_{i}))$ denotes the rank of $s(X_{i})$ in the sample $\{s(X_j),1\leq j\leq n\}$. We refer to \cite{HajSid67,vdV98} for basic
results related to linear rank statistics. The following relation is well-known:
\[
\frac{n_+n_-}{n+1}\widehat{\mbox{\sc AUC}}(s)+\frac{n_+(n_++1)}{2}=T_{n}(s)~.
\]

Moreover, the statistic $T_{n}(s)/n_+$ is an asymptotically normal estimate of
\[
W(s)=
\EXP \left(F_{s}(s(X))\mid Y=1\right)~.
\]

Note the theoretical counterpart of the previous relation may be written as
\[
W(s)=(1-p)\mbox{\sc AUC}(s)+p/2~.
\]

Now, in order to take into account a proportion $u_0$ of the highest ranks only, one may consider
the criterion related to the score generating function $\Phi_{u_0}(v)=v~\I\{v> 1-u_0\}$:
\[
W(s, u_0)
= \EXP \left(\Phi_{u_0}(F_{s}(s(X)))\mid Y=1\right)
\]
which we shall call the {\em $W$-ranking error at rate $u_0$}.

\medskip

Note that its empirical counterpart is given by $T_{n}(s, u_0)/n_+$, with
\[
T_{n}(s, u_0)=\sum_{i=1}^n\I\{Y_{i}=1\}~\Phi_{u_0}\left(\frac{\mbox{rank}(s(X_{i}))}{n+1}\right)~.
\]


Using the results from the previous subsection, we can easily check that the following theorem holds.

\begin{theorem}\label{thm:optwilcox}
We have, for all $s$:
\[
\forall s^*\in\S^*, \quad W(s, u_0) \le W(s^*, u_0)~.
\]
Furthermore, we have:
\[
W(s, u_0)=\frac{p}{2} \beta(s, u_0) (2-\beta(s, u_0)) + (1-p)\mbox{\sc LocAUC} (s, u_0)~.
\]
\end{theorem}

\begin{proof}
The result easily follows from the following representation of $\mu$:
\[
W(s, u_0) = \int_{Q(s, 1-u_0)}^{+\infty} F_s(z)~G_s(dz)
\]
and from the fact that: $F_s = p G_s + (1-p) H_s$. \qed
\end{proof}

\begin{remark}{\sc (On the choice of a score generating function $\Phi$)}
The idea of weighting the empirical AUC criterion with non-uniform
weights is equivalent to considering smooth score generating
functions $\Phi$ instead of our $\Phi_{u_0}$ in the $W$-ranking
error. Deriving optimality results for smooth criteria with our
method is straightforward but we point out that, in this case,
probabilistic interpretations are lost. In this approach, the
stochastic processes arising are rank processes for which there is
no theory available at this moment.
\end{remark}

\begin{remark}{\sc (Evidence against 'divide-and-conquer' strategies)}
It is noteworthy that not all combinations of $\beta(s,u_0)$ (or
$\alpha(s, u_0)$) and $R(s, u_0)$ lead to a criterion with $\S^*$
being the set of optimal scoring functions. We have provided two
non-trivial examples for which this is the case (Theorems
\ref{thm:optauc} and \ref{thm:optwilcox}). But, in general, this
remark should prevent from considering naive 'divide-and-conquer'
strategies for solving the local ranking problem. By naive
'divide-and-conquer' strategies, we refer here to stagewise
strategies which would, first, compute an estimate $\hat{C}$ of
the set containing the best instances, and then, solve the ranking
problem over $\hat{C}$ as described in \cite{CLV06}. However, this
idea combined with a certain amount of iterativeness might be the
key to the design of efficient algorithms. In any case, we stress
here the importance of making use of a global criterion,
synthesizing our double goal: finding and ranking the best
instances.
\end{remark}

\section{Empirical risk minimization of the local {\sc AUC} criterion}\label{sec:erm}

In the previous section, we have seen that there are various
performance measures which can be considered for the problem of
ranking the best instances. In order to perform the statistical
analysis, we will favor the representations of {\sc LocAUC} and
$W$ which involve the classification error $L(s, u_0)$ and the
local ranking error $R(s, u_0)$. By combining Theorems
\ref{thm:optauc} and \ref{thm:optwilcox}, we can easily get:
\[
2p(1-p) \mbox{\sc LocAUC}(s, u_0) = (1-p)(p+u_0)-(1-p)L(s, u_0) -
R(s, u_0)
\]
and
\[
2p W(s, u_0) = C(p, u_0) + \left(\frac{p+u_0}{2}-1\right) L(s,
u_0) - \frac{1}{4} L^2(s, u_0) - R(s, u_0)
\]
where $C(p, u_0)$ is a constant depending only on $p$ and $u_0$.

\medskip

We exploit the first expression and choose to study the
minimization of the following criterion for ranking the best
instances:
\[
M(s) \circeq M(s, u_0)=R(s,C_{s,u_0})+(1-p)L(s,u_0)~.
\]
It is obvious that the elements of $\S^*$ are the optimal elements
of the functional $M(~\cdot~, u_0)$ and we will now consider
scoring functions obtained through empirical risk minimization of
this criterion.

\medskip

More precisely, given $n$ i.i.d. copies $(X_1,Y_1), \ldots,
(X_n,Y_n)$ of $(X,Y)$, we introduce the empirical counterpart:
\[
\hat{M}_n(s) \circeq \hat{M}_n(s,
u_0)=\hat{R}_n(s)+\frac{n_-}{n}\hat{L}_n(s),
\]
with $n_-=\sum_{i=1}^n \mathbb{I}\{Y_i=-1\}$ and
\[
\hat{R}_n(s)=\frac{1}{n(n-1)}\sum_{i\neq
j}\mathbb{I}\{(s(X_i)-s(X_j))(Y_i-Y_j)>0, \; s(X_i)\wedge s(X_j)
\geq \hat{Q}(s,1-u_0) \}~.
\]

Note that $\hat{R}_n(s)$ is expected to be close to the
$U$-statistic of degree two
\[
R_n(s)=\frac{1}{n(n-1)}\sum_{i\neq j}k_s((X_i,Y_i),(X_j,Y_j)),
\]
with symmetric kernel
\[
k_s((x,y),(x',y'))=\mathbb{I}\{(s(x)-s(x'))(y-y')>0, \; s(x)\wedge
s(x') \geq Q(s,1-u_0) \}~.
\]

The statistic $R_n(s)$ corresponds to an unbiased estimate of the
local ranking error $R(s, u_0)$. The next result provides a
standard error bound for the excess risk of the empirical risk
minimizer over a  class $\mathcal{S}$ of scoring functions:
\[
\hat{s}_n=\argmin_{s\in \mathcal{S}}\hat{M}_n(s)~.
\]

\begin{proposition}
Assume that conditions (i)-(ii) of Theorem 2 are fulfilled. Then,
there exist constants $c_1$ and $c_2$ such that, for any
$\delta>0$, we have:
\[
M(\hat{s}_n)-\inf_{s\in \mathcal{S}}M(s)\leq c_1
\sqrt{\frac{V}{n}} +c_2\sqrt{\frac{\ln(1/\delta)}{n}}
\]
with probability larger than $1-\delta$.
\end{proposition}

\begin{proof}{\sc (sketch)}
The proof combines the argument used in the proof of Theorem
\ref{thm:firstorder} with the techniques used in establishing
Proposition 2 in \cite{CLV05}.
\begin{multline*}
M(\hat{s}_n)-\inf_{s\in \S}M(s) \le 2 \left(\sup_{s\in
\S}\left|\hat{R}_n(s)-R_n(s)\right| +\sup_{s\in
\S}\left|R(s)-R_n(s)\right|\right) \\
+ 2(1-p)\left( \sup_{s\in \S}\left|\hat{L}_n(s)-L_n(s)\right|+
\sup_{s\in \S}\left|L(s)-L_n(s)\right|\right) +
2\left|\frac{n_+}{n}-p\right|~.
\end{multline*}

The middle term  may be bounded by applying the result stated in
Theorem \ref{thm:firstorder}, while the last one can be handled by
using Bernstein's exponential inequality for an average of
Bernoulli random variables. By combining Lemma 1 in \cite{CLV05}
with the Chernoff method, we can deal with the $U$-process term
$\sup_{s\in \S}\left|R(s)-R_n(s)\right|$. Finally, the term
$\sup_{s\in \S}\left|\hat{R}_n(s)-R_n(s)\right|$ can also be
controlled by repeating the argument in the proof  of Theorem
\ref{thm:firstorder}. The only difference here is that we have to
consider the $U$-process term
\[
\sup_{(s,t)} \left| \frac{2}{n(n-1)}\sum_{i\neq j}
\{K_{s,t}((X_i,Y_i),(X_j,Y_j))-\mathbb{E}[K_{s,t}((X,Y),(X',Y'))]\}
\right|\,
\]
with
\[
K_{s,t}((x,y),(x',y'))=\mathbb{I}\{(s(x)-s(x'))(y-y')>0, \;
s(x)\wedge s(x') \geq t \}~.
\]
For deriving first-order results with such a process, we refer to
the same type of argument as used in \cite{CLV05}.\qed\end{proof}

\begin{remark}{\sc (about the possibility of deriving fast rates)}
By checking the proof sketch, it turns out that sharper bounds may
be achieved for the $U$-process term. Indeed, it is a simple
variation of our previous work in \cite{CLV05} where we have used
Hoeffding's decomposition in order to grasp the deep structure of
the underlying statistic. Here we will need, in addition,
condition (iii) to hold for all $u\in (0,u_0]$. Indeed, if we
localize our low-noise assumption from \cite{CLV05}, it takes the
following form: there exist constants $\alpha \in (0,1)$ and $B>0$
such that, for all $t\ge 0$, we have
\[
\forall x \in C^*_{u_0}, \qquad \Prob{\left|  \eta(X)-\eta(x)
\right| \le t} \le B \, t^{\frac{\alpha}{1-\alpha}}~.
\]
It is easy to see that this is equivalent to condition (iii) for
all $u\in (0,u_0]$: there exist constants $\alpha \in (0,1)$ and
$B>0$ such that, for all $t\ge 0$, we have
\[
\forall u \in (0,u_0], \qquad \Prob{\left|  \eta(X)-Q(\eta, 1-u)
\right| \le t} \le B \, t^{\frac{\alpha}{1-\alpha}}~.
\]
However, in the present formulation where $p$ is assumed to be
unknown, it looks like this improvement will be spoiled by the
'proportion term' which will still be of the order of a
$O(n^{-1/2})$.
\end{remark}

\section*{Appendix - Proof of Proposition \ref{prop:hoeff}}

First, for all $(s,v)\in \S\times (0,1)$ set
\[
V_n(s,v)=\frac{1}{n}\sum_{i=1}^n Y_i\I\{s(X_i)\le Q(s,
v)\}-K(s,v)~.
\]
We have the following decomposition:
\[
\forall v\in [0, 1]~, \quad  \hat{K}_n(s,v)-K(s,v)=V_n(s,F_s\circ
\hat{F}_s^{-1}(v))+ K(s,F_s\circ \hat{F}_s^{-1}(v))-K(s,v)~.
\]
We shall first prove that
\[
V_n(s,F_s\circ \hat{F}_s^{-1}(v_0))=V_n(s,v_0)+O_{\PROB}(n^{-1}).
\]
We denote by $A(s, \epsilon)$ the event $\left\{
\left|F_s\circ\hat{F}_s^{-1}(v_0)-v_0\right|  < \epsilon\right\}$.
On the event $A(s, \epsilon)$, we have:
\[
\left|V_n(s, F_s\circ\hat{F}_s^{-1}(v_0))-V_n(s,v_0)\right| \le
\sup_{v ~:~|v-v_0| <\epsilon}  \left|V_n(s, v)-V_n(s,v_0)\right|
\, .
\]
We bound the right hand side for fixed $\epsilon$, by making use
of an argument from \cite{vdG00}. First, we need to put things
into the right format. Set:
\[
V_n(s,v)-V_n(s,v_0)= \frac{1}{n}\sum_{i=1}^n
\left(u_i(s,v)-u_i(s,v_0)\right) ~,
\]
where $u_i(s,v)=Y_i \I \{s(X_i)\le Q(s, v) <0
\}-\EXP\left(Y\I\{s(X)\le Q(s, v)\}\right)$ for $s\in\S$ and $v\in
(0,1)$. We observe that
\[
\left|u_i(s,v)-u_i(s,v_0)\right|\leq d_i(v,v_0),
\]
where
\[
d_i(v,v_0) =\I \{s(X_i)\in [Q(s, v\wedge v_0), Q(s, v\vee v_0)]\}
+ |v-v_0| ~.
\]
Denote by
\[
\hat{d}(v,v_0) = \frac{1}{n}\sum_{i=1}^n \I\{s(X_i)\in [Q(s,
v\wedge v_0), Q(s, v\vee v_0)]\} + |v-v_0| ~.
\]
a distance over $\RR$. Set also:
\[
\hat{R}(\epsilon) = \sup_{v ~:~ |v-v_0|<\epsilon} \hat{d}(v,v_0)\,
.
\]
and observe that
\[
\hat{R}(\epsilon) = \frac{1}{n}\sum_{i=1}^n \I\{s(X_i)\in [Q(s,
v_0- \epsilon), Q(s, v_0+\epsilon)]\} + \epsilon \, .
\]
We then have, by applying Lemma 8.5 from \cite{vdG00}, for
$nt^2/\hat{R}^2(\epsilon)$ sufficiently large,
\[
\Prob{\sup_{v ~:~|v-v_0| \le
\epsilon}\left|V_n(s,v)-V_n(s,v_0)\right|\ge t ~\bigg|~ X_1,
\ldots, X_n }\le C \exp\left\{-
\frac{cnt^2}{\hat{R}^2(\epsilon)}\right\}~,
\]
for some positive constants $c$ and $C$. It remains to integrate
out and, for this purpose, we introduce the event:
\[
\forall x>0~, \qquad \Delta(x) = \left\{3\epsilon -x\le
\hat{R}(\epsilon) \le 3\epsilon+x \right\}.
\]
We then have:
\[
\EXP \left(\exp\left\{-
\frac{cnt^2}{\hat{R}^2(\epsilon)}\right\}\right) \le \exp\left\{-
\frac{cnt^2}{(3\epsilon+x)^2}\right\}
 + \Prob{\overline{\Delta(x)}}~.
\]
Now, we have, by Bernstein's inequality:
\[
\Prob{\overline{\Delta(x)}} = 2\Prob{\frac{1}{n} B(n, 2\epsilon) -
2\epsilon >x} \le 2 \exp \left\{ -
\frac{3nx^2}{16\epsilon}\right\}
\]
where we have used the notation $B(n, 2\epsilon)$ for a binomial
$(n, 2\epsilon)$ random variable. We can take $x =
O(t/\sqrt{\epsilon})$ and assume also $x=o(\epsilon)$ to get, for
$nt^2/\epsilon^2$ large enough,
\[
\Prob{\sup_{v ~:~|v-v_0| \le
\epsilon}\left|V_n(s,v)-V_n(s,v_0)\right|\ge t } \le C
\exp\left\{- \frac{cnt^2}{\epsilon^2}\right\}~,
\]
for some positive constants $c$ and $C$. This can be reformulated,
by writing that the following bound holds, with probability larger
than $1-\delta/2$,
\[
\sup_{v~:~|v-v_0| \le \epsilon}\left|V_n(s,v)-V_n(s,v_0)\right|\le
\epsilon \sqrt{\frac{\log(2C/\delta)}{nc}}~.
\]
We recall that, by the triangle inequality and
Dvoretsky-Kiefer-Wolfowitz theorem, if we take $\epsilon = c ~
\sqrt{\frac{\log (2/\delta)}{n}}$, we have $\Prob{A(s, \epsilon)}
\ge 1-\delta/2$. It follows that, with probability larger than
$1-\delta$, we have, for some constant $\kappa$:
\[
\left|V_n(s, F_s\circ\hat{F}_s^{-1}(v_0))-V_n(s,v_0)\right| \le
\kappa ~\left(\frac{\log(1/\delta)}{n}\right) ~,
\]
for any $s\in\S$. Now it remains to deal with the second term
$K(s,F_s\circ \hat{F}_s^{-1}(v_0))-K(s,v_0)$. Therefore, by the
differentiability assumption, we have: $\forall s \in \S$,
\[
\sup_{|v-v_0| \le
\delta}\{K(s,v)-K(s,v_0)-(v-v_0)K'(s,v_0)\}=O(\delta^2)~, \quad
\mbox{as} ~ \delta \rightarrow 0~.
\]
Since $\vert F_s\circ
\hat{F}_s^{-1}(v_0))-v_0\vert=O_{\PROB}(n^{-1/2})$, we get that
\[
K(s,F_s\circ \hat{F}_s^{-1}(v_0))-K(s,v_0)=K'(s,v_0)(F_s\circ
\hat{F}_s^{-1}(v_0)-v_0)+O_{\mathbb{P}}(n^{-1})~, \quad \mbox{as}
~ n\rightarrow \infty~.
\]
Moreover, as
\[
F_s\circ \hat{F}_s^{-1}(v_0)-v_0=-(\hat{F}_s\circ
F_s^{-1}(v_0)-v_0)+O_{\mathbb{P}}(n^{-1})~,
\]
we finally obtain that
\[
K(s,F_s\circ
\hat{F}_s^{-1}(v_0))-K(s,v_0)=-K'(s,v_0)(\hat{F}_s\circ
F_s^{-1}(v_0)-v_0)+O_{\mathbb{P}}(n^{-1})~. \mbox{\qed}
\]

\bigskip

\noindent {\bf Acknowledgements. } We thank St\'ephane Boucheron
and G\'abor Lugosi for their helpful remarks and encouragements.

\bibliography{References_ranking}
\bibliographystyle{plain}

\end{document}